\documentclass[12pt]{article}
\usepackage{amssymb,amsmath}

\def\ln{\mbox{ln\,}}

\def\stac{\stackrel}
\def\~{\stac{\sim}}
\def\8{\infty}

\def\({\left(}
\def\){\right)}

\def\be{\begin{equation}}
\def\ee{\end{equation}}

\def\ba{\begin{array}}
\def\ea{\end{array}}

\def\({\left(}\def\){\right)}

\begin{document}
\title{
Quenching, global existence and blowup phenomena in heat transfer
}
   \author{
Xianfa Song{\thanks{E-mail:\ \tt songxianfa2004@163.com
  }}\\
 \small Department of Mathematics, School of Mathematics, Tianjin University,\\
\small Tianjin, 300072, P. R. China
}

\maketitle
\date{}

\newtheorem{theorem}{Theorem}[section]
\newtheorem{definition}{Definition}[section]
\newtheorem{lemma}{Lemma}[section]
\newtheorem{proposition}{Proposition}[section]
\newtheorem{corollary}{Corollary}[section]
\newtheorem{remark}{Remark}[section]
\renewcommand{\theequation}{\thesection.\arabic{equation}}
\catcode`@=11 \@addtoreset{equation}{section} \catcode`@=12

\begin{abstract}

In this paper, we are concerned with the following problem appearing in heat transfer:
\begin{equation*}
\left\{
\begin{array}{llll}
&\frac{\partial u}{\partial t}=d_1\Delta u-a(x)\cdot \nabla u+f(u,v),\quad x\in\Omega,\ t>0,\\
&\frac{\partial v}{\partial t}=d_2\Delta v-b(x) \cdot \nabla v+g(u,v),\quad x\in\Omega,\ t>0,\\
&\frac{\partial u}{\partial \eta}=\frac{\partial v}{\partial \eta}=0\quad {\rm or} \quad u=v=0, \qquad x\in\partial\Omega,\ t>0,\\
&u(x,0)=u_0(x),\quad v(x,0)=v_0(x),\quad x\in\Omega.
\end{array}\right.
\end{equation*}
Here $\Omega\subset \mathbb{R}^N(N\geq 1)$, is a bounded smooth domain, $d_1, d_2>0$, $f(u,v)$ and $g(u,v)$ are smooth functions of $(u,v)$, $u_0(x)$ and $v_0(x)$ are nonnegative continuous functions of $x$, while $a(x)$, $b(x)$ are vector valued functions.  Basing on the relations between a system of ODE and a system of parabolic equations, we establish some general theories in heat transfer about quenching, global existence and blowup phenomena, obtain the conditions(even watershed) on $f(u,v)$, $g(u,v)$, $a(x)$ and $b(x)$ which let the solution be global existence, quench or blow up, and estimate the bounds for blowup time and quenching time.

{\bf Keywords:} Heat transfer; Sub-controlling(or sup-controlling) method; Quenching; Global existence; Blow up.

{\bf 2010 MSC:} 35K55; 35K58.

\end{abstract}

\section{Introduction}

\qquad In this paper, we are concerned with the following problem appearing in heat transfer:
\begin{equation}
\label{1}\left\{
\begin{array}{llll}
&\frac{\partial u}{\partial t}=d_1\Delta u-a(x)\cdot \nabla u+f(u,v),\quad x\in\Omega,\ t>0,\\
&\frac{\partial v}{\partial t}=d_2\Delta v-b(x) \cdot \nabla v+g(u,v),\quad x\in\Omega,\ t>0,\\
&\frac{\partial u}{\partial \eta}=\frac{\partial v}{\partial \eta}=0\quad {\rm or}\quad u=v=0,  \qquad x\in\partial\Omega,\ t>0,\\
&u(x,0)=u_0(x),\quad v(x,0)=v_0(x),\quad x\in\Omega.
\end{array}\right.
\end{equation}
Here $\Omega\subset \mathbb{R}^N(N\geq 1)$, is a bounded smooth domain,  $d_1, d_2>0$, $f(u,v)$ and $g(u,v)$ are smooth functions of $(u,v)$, $u_0(x)$ and $v_0(x)$ are nonnegative continuous functions of $x$, while $a(x)$, $b(x)$ are vector valued functions. Model (\ref{1}) often appears in heat transfer, $u$ and $v$ represent the temperature, the terms $a(x)\cdot \nabla u$ and $b(x)\cdot \nabla v$ are called convection terms, which means that there exists convection in the course of heat transfer. The local wellposedness of a parabolic system such as (\ref{1}) was proved under certain assumptions on $f(u,v)$, $g(u,v)$, $a(x)$, $b(x)$, $u_0(x)$ and $v_0(x)$(see \cite{Evans} and the references therein). In convenience, we denote the problem (\ref{1}) subject to Neumann boundary condition by (\ref{1}A) and (\ref{1}) subject to Dirichlet boundary condition by (\ref{1}B).

The motivations of this paper are as follows.

First, we often take a system of ODE as a sub-controlling or sup-controlling system of parabolic equations, and call the solution of ODE as the sub-solution or sup-solution of the system of parabolic equations. About the recent results on the controllability of parabolic equations, we can refer to \cite{Doubova, Duyckaerts, Fernandez, Privat, Zuazua} and the references therein.  Consider the following ODE problem:
\begin{equation}
\label{1'}\left\{
\begin{array}{llll}
&\frac{dw}{dt}=f(w,z),\quad \frac{dz}{dt}=g(w,z),\quad t>0,\\
&w(0)=c_1\geq 0,\quad z(0)=c_2\geq 0.
\end{array}\right.
\end{equation}
We try to reveal the relations between (\ref{1'}) and (\ref{1}) when they contain the same functions $f$ and $g$. Roughly speaking, if the solution of (\ref{1'}) is global existence(or blows up in finite time, or quenches in finite time) under certain assumptions on $f$ and $g$, then we hope to prove that the solution of (\ref{1}) is also global existence(or blows up in finite time, or quenches in finite time) under the same assumptions on $f$ and $g$ with suitable conditions on $a(x)$ and $b(x)$, which is a basic and very interesting question in this direction.

{\bf Remark 1.1.} In this paper, (\ref{1}) and (\ref{1'}) contain the same functions $f$ and $g$ means that the expressions of $f$ and $g$ in (\ref{1}) and those in (\ref{1'}) are the same. For example, if $f(u,v)=u^{p_1}v^{q_1}$  and $g(u,v)=u^{p_2}v^{q_2}$ in (\ref{1}), then (\ref{1'}) contains the same functions $f$ and $g$ means that $f(w,z)=w^{p_1}z^{q_1}$ and $g(w,z)=w^{p_2}z^{q_2}$.

Second, many authors studied the following problem:
\begin{equation}
\label{h2}\left\{
\begin{array}{llll}
&\frac{\partial u}{\partial t}=d_1\Delta u+f(u,v),\quad \frac{\partial v}{\partial t}=d_2\Delta v+g(u,v),\quad x\in \Omega,\quad t>0,\\
&\frac{\partial u}{\partial \eta}=\frac{\partial v}{\partial \eta}=0\quad {\rm or}\quad u=v=0,\qquad x\in \partial \Omega,\ t>0,\\
&u(x,0)=u_0(x),\quad v(x,0)=v_0(x),\quad x\in \Omega.
\end{array}\right.
\end{equation}
For some special $f$ and $g$ such as $f(u,v)=a_1v^{p_1}u^{q_1}$, $g(u,v)=a_2v^{p_2}u^{q_2}$ or $f(u,v)=F(u)\exp(\frac{v-1}{v})$, $g(u,v)=F(u)\exp(\frac{v-1}{v})$, they obtain the conditions which can ensure that the solution of (\ref{h2}) is global existence or blows up in finite time or quenches in finite time. We can refer to \cite{Bandle, Bebernes, Castillo, Escobedo, Fujishima, Gang, Hu, Ishige, Li, Pierre, Rossi, Wang} and the references therein. However, there are few results on (\ref{h2}) about such conditions on more general $f$ and $g$. Naturally, we are concerned with the following question: What conditions on general $f$ and $g$ can make the solutions of (\ref{1}), (\ref{1'}) and (\ref{h2}) be global existence, or blow up in finite time, or quench in finite time? Meanwhile, we are interested in the following question: How the convection terms $a(x)\cdot \nabla u$ and $b(x)\cdot \nabla v$  effect the properties for the solution of (\ref{1})?

Basing on the motivations above, we hope to establish the conditions on general $f$ and $g$ which can make the solutions of (\ref{1}), (\ref{1'}) and (\ref{h2}) be global existence, or blow up in finite time, or quench in finite time.

We would like to say something about the lectures on blowup phenomenon of parabolic equations. Blowup for nonlinear evolution equations
had been deserved a great deal of interest ever since the pioneering papers \cite{Fujita1, Fujita2}. For the following scalar semi-linear parabolic equation problem
\begin{equation}
\label{h1}\left\{
\begin{array}{llll}
&u_t=d\Delta u+\hat{f}(u),\quad x\in \Omega,\quad t>0,\\
&\frac{\partial u}{\partial \eta}=0\quad {\rm or}\ u=0,\quad x\in \partial \Omega,\ t>0,\\
&u(x,0)=u_0(x),\quad x\in \Omega,
\end{array}\right.
\end{equation}
the key condition which makes the solution blow up in finite time is $$\int^{+\infty}_c\frac{du}{\hat{f}(u)}<+\infty,\quad c>0.$$ We can refer to \cite{Ball, Bebernes1, Bellout, Escher, Lacey, Levine1, Meier, Weissler} and the references therein.  Some authors also considered the problems of scalar parabolic equation subject to $\frac{\partial u}{\partial \eta}=\hat{f}(u)$, or a system of parabolic equations subject to nonlinear boundary conditions $\frac{\partial u}{\partial \eta}=r(u,v)$ and $\frac{\partial v}{\partial \eta}=s(u,v)$. We can refer to \cite{Acosta, Gomez, Rial} and the references therein.

About the quenching phenomenon for the solution of a parabolic equation, it was first discussed by Kawarada in \cite{Kawarada}. We also can refer to \cite{Dai1, Kalashnik, Levine2, Levine3} and the references therein to see more information.

Since the complexity of interaction between different nonlinearities, for a system of nonlinear parabolic equations, non-simultaneous blowup or non-simultaneous quenching phenomenon may happen, we can refer to \cite{Pablo, Quiros1, Quiros2, Zheng} and the references therein. However, we don't establish general results on the non-simultaneous blowup or non-simultaneous quenching phenomenon in the framework of theory here, we only give some examples to show the non-simultaneous blowup or non-simultaneous quenching phenomenon.

We don't state the precise expressions of our results here in order to control the length of this section, we will state and prove them in the corresponding sections. But we would like to compare our results with others below.

1. (i). Our results on blowup and global existence of the solutions to (\ref{1'}) and (\ref{h2}) meet with those of others. However, we establish some basic and interesting results on (\ref{1}), (\ref{1'}) and (\ref{h2}) for general functions $f$ and $g$ in the framework of theory and find more interesting phenomena. For instance, from Example 3.2, we see that: some nonlinearities $f(u,v)$ and $g(u,v)$ make $(u,v)$ blow up but $(u^{\alpha}v^{\beta})$ exist globally, while some nonlinearities $f(u,v)$ and $g(u,v)$ can make both $(u,v)$ and $(u^{\alpha}v^{\beta})$ blow up.

(ii) There are few results on (\ref{SM11'}) and (\ref{SM11}) before. However, by the results of Theorem A and Theorem 1 in this paper, whether the solution is global existence or blows up in finite time is determined by $m\leq n$ or $m>n$, the watershed is $m=n$.

2. Differing to the nonlinear Neumann boundary conditions $\frac{\partial u}{\partial \eta}=r(u,v)$ and $\frac{\partial v}{\partial \eta}=s(u,v)$ under the assumptions
of $r(u,v)>0$ and $s(u,v)>0$, and the properties for the solution are very dependent on $r(u,v)$ and $s(u,v)$ in \cite{Acosta}, (\ref{1}) and (\ref{h2}) have homogeneous Neumann or Dirichlet boundary conditions in this paper.

3. We establish the results on quenching phenomena for the solutions of (\ref{1}) and (\ref{h2}) with general $f(u,v)$ and $g(u,v)$, while others considered (\ref{h2}) with some special $f(u,v)$ and $g(u,v)$ such as $f(u,v)=u^{-p_1}v^{-q_1}$ and $g(u,v)=u^{-p_2}v^{-q_2}$.

4. We consider the roles of convection terms $a(x)\cdot \nabla u$ and $b(x)\cdot \nabla v$, and find that some type of convection terms can delay the quenching time, while some type of convection terms can make the solution exist globally.

The rest of this paper is organized as follows. In Section 2, we will give the results on (\ref{1}), (\ref{1'}) and (\ref{h2}) when $f$ and $g$ are separation of variables. In Section 3, we will give the results on (\ref{1}), (\ref{1'}) and (\ref{h2}) when $f$ and $g$ aren't separation of variables. In Section 4, we will study the quenching phenomenon. In Section 5, we will consider the roles of the convection terms $a(x)\cdot \nabla u$ and $b(x)\cdot \nabla v$.

\section{The results on (\ref{1}) when $f$ and $g$ are separation of variables}
\qquad In this section, we deal with (\ref{1}) in the special case of $f(u,v)=f_1(v)g_1(u)$ and $g(u,v)=f_2(u)g_2(v)$, i.e., $f$ and $g$ are separation of variables. We hope to judge whether the solution is global existence or blowup in finite time directly by the structures of $f$ and $g$.

\subsection{Theorem A and Theorem 1}
\qquad Consider the following ODE problem:
\begin{equation}
\label{661}\left\{
\begin{array}{llll}
&w_t=f_1(z)g_1(w),\quad z_t=f_2(z)g_2(w),\quad t>0,\\
&w(0)=c_1\geq 0,\quad z(0)=c_2\geq 0.
\end{array}\right.
\end{equation}
Here $f_i(z)\geq 0(\not\equiv 0)$, $g_i(w)\geq 0(\not\equiv 0)$, $i=1,2$,  $ f'_1(z)\geq 0(\not\equiv 0)$ and $g'_2(w)\geq 0(\not\equiv 0)$.

Assume that $\int^{+\infty}_{c_1}\frac{dw}{g_1(w)}=+\infty$ and $\int^{+\infty}_{c_2}\frac{dz}{f_2(z)}=+\infty$. Let
\begin{align}
\int^{w(t)}_{c_1}\frac{ds}{g_1(s)}=G_1(w(t)),\quad \int^{z(t)}_{c_2}\frac{d\theta}{f_2(\theta)}=F_2(z(t)).\label{662}
\end{align}
Then (\ref{661}) becomes
\begin{equation}
\label{9121'}\left\{
\begin{array}{llll}
&(G_1(w))_t=f_1(z),\quad (F_2(z))_t=g_2(w),\quad t>0,\\
&w(0)=c_1\geq 0,\quad z(0)=c_2\geq 0.
\end{array}\right.
\end{equation}
Let
\begin{align}
&W(t)=G_1(w(t)),\quad Z(t)=F_2(z(t)),\label{663}\\
&\tilde{f}(Z)=f_1[F^{-1}_2(Z)],\quad \tilde{g}(W)=g_2[G^{-1}_1(W)].\label{9124}
\end{align}
(\ref{661}) becomes
\begin{equation}
\label{661''}\left\{
\begin{array}{llll}
&W_t=\tilde{f}(Z),\quad Z_t=\tilde{g}(W),\quad t>0,\\
&W(0)=\tilde{c}_1\geq 0,\quad Z(0)=\tilde{c}_2\geq 0.
\end{array}\right.
\end{equation}

We have

{\bf Theorem A.} {\it Assume that $f(w,z)=f_1(z)g_1(w)$ and $g(w,z)=f_2(z)g_2(w)$, $f_i(z), f'_i(z)\geq 0$ for $z\geq 0$ and $f_i(z), f'_i(z)\not\equiv 0$
in any subinterval of $(0,+\infty)$, $g_i(w), g'_i(w)\geq 0$ for $w\geq 0$ and $g_i(w), g'_i(w)\not\equiv 0$
in any subinterval of $(0,+\infty)$, $i=1,2$.

(1). If $\int^{+\infty}_{c_1}\frac{dw}{g_1(w)}<+\infty$ or $\int^{+\infty}_{c_2}\frac{dz}{f_2(z)}<+\infty$. Then the solution of (\ref{661})
will blow up in finite time for positive initial data.

(2). Assume that $\int^{+\infty}_{c_1}\frac{dw}{g_1(w)}=+\infty$ and $\int^{+\infty}_{c_2}\frac{dz}{f_2(z)}=+\infty$. And there exist $0\leq \tilde{c}_3\leq \tilde{c}_1$, $0\leq \tilde{c}_4\leq \tilde{c}_2$,
$\tilde{f}(Z), \tilde{f}'(Z)\geq 0$ for all $Z\geq \tilde{c}_3$, $\tilde{g}(W), \tilde{g}'(W)\geq 0$ for all $W\geq \tilde{c}_4$. Let $\tilde{F}(Z(t))=\int^{Z(t)}_{\tilde{c}_4}\tilde{f}(s)ds$ and $\tilde{G}(W(t))=\int^{W(t)}_{\tilde{c}_3}\tilde{g}(\theta)d\theta$ for $t\geq 0$. Suppose
that there exist positive constants $\epsilon$ and $K$ such that
\begin{align}
\epsilon \tilde{G}(\tilde{c}_1)\leq \tilde{F}(\tilde{c}_2)\leq K\tilde{G}(\tilde{c}_1).\label{12311'}
\end{align}
Then the solution of (\ref{661''}) is global existence if
\begin{align}
\int_{\tilde{c}_1}^{+\infty}\frac{ds}{\tilde{f}[\tilde{F}^{-1}(K\tilde{G}(s))]}=+\infty\quad{\rm and}\quad  \int_{\tilde{c}_2}^{+\infty}\frac{d\theta}{\tilde{g}[\tilde{G}^{-1}(\frac{1}{\epsilon}\tilde{F}(\theta))]}=+\infty,\label{19111'}
\end{align}
and it will blow up in finite time for large initial data if
\begin{align}
\int_{\tilde{c}_1}^{+\infty}\frac{ds}{\tilde{f}[\tilde{F}^{-1}(\epsilon\tilde{G}(s))]}<+\infty\quad{\rm or}\quad  \int_{\tilde{c}_2}^{+\infty}\frac{d\theta}{\tilde{g}[\tilde{G}^{-1}(\frac{1}{K}\tilde{F}(\theta))]}<+\infty.\label{19112'}
\end{align}
Here $\tilde{F}^{-1}$ and $\tilde{G}^{-1}$ are the inverse functions of $\tilde{F}$ and $\tilde{G}$ respectively. Therefore, the solution of (\ref{661})
is global existence or blows up in finite time in the corresponding case.

}

Parallel to Theorem A, we have the following conclusions on (\ref{1}) and (\ref{h2}).

{\bf Theorem 1.} {\it  Assume that $f(u,v)=f_1(v)g_1(u)$ and $g(u,v)=f_2(v)g_2(u)$, $f_i(v), f'_i(v)\geq 0$ for $v\geq 0$ and $f_i(v), f'_i(v)\not\equiv 0$
in any subinterval of $(0,+\infty)$, $g_i(u), g'_i(u)\geq 0$ for $u\geq 0$ and $g_i(u), g'_i(u)\not\equiv 0$
in any subinterval of $(0,+\infty)$, $i=1,2$.

(1). If the assumptions of Theorem A(2) and (\ref{19111'}) hold, then the solutions of (\ref{1}A),(\ref{1}B), (\ref{h2}A) and (\ref{h2}B) are global existence
for any nonnegative initial data $(u_0,v_0)$.

(2). Suppose that
\begin{align}
&d_1\Delta u_0-a(x)\cdot \nabla u_0+f(u_0,v_0)\geq 0\quad {\rm for}\quad  x\in\Omega,\label{12041}\\
&d_2\Delta v_0-b(x)\cdot \nabla v_0+g(u_0,v_0)\geq 0\quad {\rm for}\quad  x\in\Omega.\label{12042}
\end{align}

If the assumptions of Theorem A(1) hold or the assumptions of Theorem A(2) hold, and (\ref{19112'}) are true,
 then the solutions of (\ref{1}A) and (\ref{h2}A) will blow up in finite time for initial data $(u_0,v_0)\geq (c_1, c_2)$. If the solution of (\ref{661}) blows up in finite time for initial data $(c_1,c_2)=\mathbf{0}$, then the solutions of (\ref{1}B) and (\ref{h2}B)  will blow up in finite time for any nonnegative initial data $(u_0,v_0)$. Here $(x_1,y_1)\geq (x_2,y_2)$ means that $x_1\geq x_2$ and $y_1\geq y_2$.
}

\subsection{A special case of Theorem A and some related results}
\qquad In order to prove Theorem A and Theorem 1, we first consider a special case of (\ref{1}),
\begin{equation}
\label{F-G}\left\{
\begin{array}{llll}
&u_t=d_1\Delta u-a(x)\cdot \nabla u+f(v),\quad x\in\Omega,\ t>0,\\
&v_t=d_2\Delta v-b(x)\cdot \nabla v+g(u),\quad x\in\Omega,\ t>0,\\
&\frac{\partial u}{\partial \eta}=\frac{\partial v}{\partial \eta}=0\quad {\rm or}\quad u=v=0,  \quad x\in\partial\Omega,\ t>0,\\
&u(x,0)=u_0(x),\quad v(x,0)=v_0(x),\quad x\in\Omega.
\end{array}\right.
\end{equation}
Here $f(v), f'(v)\geq 0$ for $v\geq 0$, $g(u), g'(u)\geq 0$ for $u\geq 0$. In convenience, we also denote the problem (\ref{F-G}) subject to Neumann boundary condition by (\ref{F-G}A) and (\ref{F-G}) subject to Dirichlet boundary condition by (\ref{F-G}B). We hope to establish the conditions on the blowup in finite time and global existence of the solution to (\ref{F-G}), which can be judged by the structures of $f(v)$ and $g(u)$ directly. To do this, we consider the following ODE problem:
\begin{equation}
\label{FG}\left\{
\begin{array}{llll}
&w_t=f(z),\quad z_t=g(w),\quad t>0,\\
&w(0)=c_1\geq 0,\quad z(0)=c_2\geq 0.
\end{array}\right.
\end{equation}
We have

{\bf Theorem A$'$} {\it 1. Assume that there exist $0\leq c_3\leq c_1$, $0\leq c_4\leq c_2$ such that
$f(z),\ f'(z)\geq 0$ for $z\geq c_3$, and $f(z),\ f'(z)\not\equiv 0$ in any subinterval of $[0,+\infty)$, $g(w),\ g'(w)\geq 0$ for $w\geq c_4$, $g(w),\ g'(w)\not\equiv 0$ in any subinterval of $[0,+\infty)$. Let
$$F(z(t))=\int^{z(t)}_{c_4}f(s)ds,\quad G(w(t))=\int^{w(t)}_{c_3}g(\theta)d\theta,\quad t\geq 0.$$
Suppose that
there exist positive constants $0<\epsilon<1$ and $K>1$ such that
\begin{align}
\epsilon G(c_1)\leq F(c_2)\leq KG(c_1).\label{12311}
\end{align}
Then the solution of (\ref{FG}) is global existence for any initial data $(c_1,c_2)$ if
\begin{align}
\int_{c_1}^{+\infty}\frac{ds}{f[F^{-1}(KG(s))]}=+\infty\quad{\rm and}\quad  \int_{c_2}^{+\infty}\frac{d\theta}{g[G^{-1}(\frac{1}{\epsilon}F(\theta))]}=+\infty,\label{19111}
\end{align}
and it will blow up in finite time for large initial data $(c_1,c_2)$ if
\begin{align}
\int_{c_1}^{+\infty}\frac{ds}{f[F^{-1}(\epsilon G(s))]}<+\infty\quad{\rm or}\quad  \int_{c_2}^{+\infty}\frac{d\theta}{g[G^{-1}(\frac{1}{K}F(\theta))]}<+\infty.\label{19112}
\end{align}
Here $F^{-1}$ and $G^{-1}$ are the inverse functions of $F$ and $G$ respectively.

2. (1). Under the condition (\ref{19111}), then the solutions of (\ref{F-G}A) and (\ref{F-G}B) are global existence for any nonnegative initial data $(u_0,v_0)$;

(2). Suppose that
\begin{align}
&d_1\Delta u_0-a(x)\cdot \nabla u_0+f(v_0)\geq 0\quad {\rm for}\quad  x\in\Omega,\label{12043}\\
&d_2\Delta v_0-b(x)\cdot \nabla v_0+g(u_0)\geq 0\quad {\rm for}\quad  x\in\Omega.\label{12044}
\end{align}

Under the condition (\ref{19112}), if the solution of (\ref{FG}) blows up in finite time for initial data $(c_1,c_2)\geq \mathbf{0}$, then the solution of (\ref{F-G}A) will blow up in finite time for initial data $(u_0,v_0)\geq (c_1, c_2)$; Parallelly, if the solution of (\ref{FG}) blows up in finite time for initial data $(c_1,c_2)=\mathbf{0}$, then the solution of (\ref{F-G}B) will blow up in finite time for any nonnegative initial data $(u_0,v_0)$.

}

{\bf Proof:} 1. First, noticing that $f(z)\geq 0$ and $g(w)\geq 0$, we know that $F(z)$, $G(w)$, $F^{-1}$ and $G^{-1}$ are all increasing functions.

Let $$
I(t)=F(z(t))-\epsilon G(w(t)),\quad J(t)=F(z(t))-KG(w(t)).
$$
Then $I(0)\geq 0$, $J(0)\leq 0$ and
$$
I'(t)=(1-\epsilon)f(z(t))g(w(t))\geq 0,\quad J'(t)=(1-K)f(z(t))g(w(t))\leq 0.
$$
Consequently,
$$
\epsilon G(w(t))\leq F(z(t))\leq KG(w(t)),\quad \frac{1}{K} F(z(t))\leq G(w(t))\leq \frac{1}{\epsilon} F(z(t))
$$
and
\begin{align*}
F^{-1}\left(\epsilon G(w(t))\right)\leq z(t)\leq F^{-1}\left(KG(w(t))\right),\\
 G^{-1}\left(\frac{1}{K} F(z(t))\right)\leq w(t)\leq G^{-1}\left(\frac{1}{\epsilon} F(z(t))\right).
\end{align*}
$f'(z)\geq 0$ and $g'(w)\geq 0$( and $\not\equiv 0$ in any subinterval of $(0,+\infty)$) imply that $f(z)$ and $g(w)$ are increasing functions, so
\begin{align*}
f[F^{-1}\left(\epsilon G(w(t))\right)]\leq f(z(t))\leq f[F^{-1}\left(K G(w(t))\right)],\\
g[G^{-1}\left(\frac{1}{K} F(z(t))\right)]\leq g(w(t))\leq g[G^{-1}\left(\frac{1}{\epsilon} F(z(t))\right)].
\end{align*}
That is,
\begin{align}
f[F^{-1}\left(\epsilon G(w(t))\right)]\leq w_t\leq f[F^{-1}\left(K G(w(t))\right)],\label{19113}\\
g[G^{-1}\left(\frac{1}{K} F(z(t))\right)]\leq z_t\leq g[G^{-1}\left(\frac{1}{\epsilon} F(z(t))\right)].\label{19114}
\end{align}
Using (\ref{19113}) and (\ref{19114}), by the sub-solution and sup-solution theory of ODE, it is easy to verify that
the solution of (\ref{1'}) is global existence if
\begin{align*}
\int_c^{+\infty}\frac{ds}{f[F^{-1}(KG(s))]}=+\infty\quad{\rm and}\quad  \int_c^{+\infty}\frac{d\theta}{g[G^{-1}(\frac{1}{\epsilon}F(\theta))]}=+\infty,
\end{align*}
and it will blow up in finite time for large initial data if
\begin{align*}
\int_c^{+\infty}\frac{ds}{f[F^{-1}(\epsilon G(s))]}<+\infty\quad{\rm or}\quad  \int_c^{+\infty}\frac{d\theta}{f[G^{-1}(\frac{1}{K}F(\theta))]}<+\infty.
\end{align*}

2. (1). Under the conditions of (\ref{19111}), the solution of (\ref{1'}) with $c_1=\max_{x\in \bar{\Omega}}u_0(x)$ and $c_2=\max_{x\in \bar{\Omega}}v_0(x)$
can be taken as a sup-solution of (\ref{F-G}), which implies that the solution of (\ref{F-G}) is global existence for any nonnegative initial data $(u_0,v_0)$.

(2). Suppose that (\ref{12043}) and (\ref{12044}) hold. Under the conditions of (\ref{19112}), if the solution of (\ref{1'}) will blow up in finite time for initial data $(c_1,c_2)$, then the solution of (\ref{1'}) with $c_1\leq c'_1=\min_{x\in \bar{\Omega}}u_0(x)$ and $c_2\leq c'_2=\min_{x\in \bar{\Omega}}v_0(x)$
can be taken as a sub-solution of (\ref{F-G}A), which implies that the solution of (\ref{F-G}A)  will blow up in finite time for initial data $u_0(x)\geq c_1$ and $v_0(x)\geq c_2$.

If the solution of (\ref{1'}) with initial data $(c_1,c_2)=\mathbf{0}$ will blow up in finite time, then we can take the sub-solution in the form of
$$(\underline{u}(x,t),\ \underline{v}(x,t))=(\underline{u}(\alpha t \delta(x)),\ \underline{v}(\beta t \xi(x))).$$ Here the functions $\underline{u}(\cdot)$ and $\underline{v}(\cdot)$ satisfy
\begin{align*}
&\frac{d\underline{u}}{dt}=f(\underline{v}),\quad \frac{d\underline{v}}{dt}=g(\underline{u}),\quad t>0,\\
&\underline{u}(0)=0,\quad \underline{v}(0)=0,
\end{align*}
i.e., $(\underline{u}(\cdot),\underline{v}(\cdot))$ is the blowup solution of (\ref{1'}) with initial data $(c_1,c_2)=\mathbf{0}$. And $\alpha$, $\beta$ are small positive constants to be determined later, the functions $\delta(x)$ and $\xi(x)$ are arbitrary independent nonnegative functions satisfy
\begin{align}
d_1\Delta \delta-a(x)\cdot \nabla \delta+p(\delta)=0\quad {\rm in}\  \Omega,\quad \delta=0\quad {\rm on}\ \partial \Omega \label{911x1}
\end{align}
and
\begin{align}
d_2\Delta \xi-b(x)\cdot \nabla \xi+q(\xi)=0\quad {\rm in}\  \Omega,\quad \xi=0\quad {\rm on}\ \partial \Omega,\label{911x2}
\end{align}
where $p(\delta)>0$ and $q(\xi)>0$ are continuous functions such that (\ref{911x1}) and (\ref{911x2}) have nonnegative solutions. Obviously, $\underline{u}(x,t)=\underline{v}(x,t)=0$ for $x\in \partial \Omega$, $t>0$. After some computations,
we have
\begin{align}
&\quad\underline{u}_t-d_1\Delta \underline{u}+a(x)\cdot \nabla \underline{u}\nonumber\\
&=\alpha f(\underline{v}(t\beta \xi(x)))[\delta(x)+tp(\delta)]-t^2d_1\alpha\beta g(\underline{u})f'(\underline{v})\nabla \delta\cdot \nabla \xi\nonumber\\
&\leq f(\underline{v}(t\beta \xi(x)))=f(\underline{v}),\label{911w1}\\
&\quad\underline{v}_t-d_2\Delta \underline{v}+b(x)\cdot \nabla \underline{v}\nonumber\\
&=\beta g(\underline{u}(t\alpha \delta(x)))[\xi(x)+tq(\xi)]-t^2d_2\alpha \beta f(\underline{v})g'(\underline{u})\nabla \delta\cdot \nabla \xi\nonumber\\
&\leq g(\underline{u}(t\alpha \delta(x)))=g(\underline{u})\label{911w2}
\end{align}
for $x\in \Omega$ and $t$ small enough if $\alpha<<1$, $\beta<<1$. That is, $(\underline{u}, \underline{v})$ is a sub-solution of (\ref{1}B).
Therefore, there exist some $x_0\in \Omega$ and $T>0$ such that
$$
\lim_{t\rightarrow T^-} [\underline{u}(x_0,t)+\underline{v}(x_0,t)]=+\infty
$$
and
$$
T\leq \frac{T^*}{\max(\alpha \delta_0, \beta \xi_0)},
$$
where $T^*$ is the blowup time of the solution to (\ref{1'}), $\delta_0=\max_{x\in \Omega}\delta(x)$ and $\xi_0=\max_{x\in \Omega}\xi(x)$.\hfill $\Box$

We would like to give some examples to illustrate the results of Theorem A$'$.

{\bf Example 2.1.} $f(z)=z^p$, $g(w)=w^q$. Let $c_3=c_4=0$. Then $F(z)=\frac{z^{p+1}}{p+1}$, $G(w)=\frac{w^{q+1}}{q+1}$, $$f[F^{-1}(KG(s))]=C_1(p,q,K)s^{\frac{p(q+1)}{p+1}},\quad
f[F^{-1}(\epsilon G(s))]=C'_1(p,q,\epsilon)s^{\frac{p(q+1)}{p+1}},$$
$$g[G^{-1}(\frac{1}{\epsilon}F(\theta))]=C'(p,q,\epsilon)\theta^{\frac{q(p+1)}{q+1}},\quad g[G^{-1}(\frac{1}{K}F(\theta))]=C'_2(p,q,K)\theta^{\frac{q(p+1)}{q+1}}.$$
The solutions of (\ref{FG}) is global existence if and only if $pq\leq 1$. Consequently, the solution of (\ref{F-G}) is global existence for any initial data if $pq\leq 1$, while the solution of (\ref{F-G}A) will blow up in finite time if $pq>1$
for $w_0(x)\geq c_1>0$ and $z_0(x)\geq c_2>0$.

{\bf Example 2.2.} $f(z)=e^z$, $g(w)=e^w$. Let $c_3=c_4=0$. Then $F(z)=e^z-1\geq \frac{e^z}{2}$ for $z>\ln 2$, $G(w)=e^w-1\geq \frac{e^w}{2}$ for $w>\ln 2$,
$$f[F^{-1}(KG(s))]\geq C_1(K)e^s,\quad
f[F^{-1}(\epsilon G(s))]\geq C'_1(\epsilon)e^s,$$
$$g[G^{-1}(\frac{1}{\epsilon}F(\theta))]\geq C_2(\epsilon)e^{\theta},\quad g[G^{-1}(\frac{1}{K}F(\theta))]\geq C'_2(K)e^{\theta}$$
for $s$, $\theta$ large enough. The solutions of (\ref{FG})  and (\ref{F-G}) will blow up in finite time
for nonnegative initial data.

{\bf Example 2.3.} (i) $f(z)=e^z$, $g(w)=\ln w$. Let $c_3=0$, $c_4=1$. Then $F(z)\leq e^z$, $G(w)=w\ln w+1-w$,
$$f[F^{-1}(KG(s))]\leq C_1(K)[s\ln s+1-s],\quad
f[F^{-1}(\epsilon G(s))]\leq C'_1(\epsilon)[s\ln s+1-s].$$
Obviously,
$$
\int^{+\infty}_1\frac{ds}{f[F^{-1}(KG(s))]}\geq\int^{+\infty}_1\frac{ds}{C_1(K)[s\ln s+1-s]}=+\infty.
$$
Although we cannot give the explicit expression of $G^{-1}$, since $G(w)\leq w\ln w$ for $w>1$, we can obtain
$$g[G^{-1}(\frac{1}{\epsilon}F(\theta))]\leq C'_1(p,q,\epsilon)\theta (1+\ln \theta)\quad {\rm for}\quad \theta>>1.$$ And
$$
\int^{+\infty}_C\frac{d\theta}{g[G^{-1}(\frac{1}{\epsilon}F(\theta))]}=+\infty.
$$
So the solution of (\ref{FG})  is global existence for initial data $c_1>1$ and $c_2\geq 0$. In fact, we even can construct the following sup-solution of (\ref{FG})
$$
\bar{w}(t)=e^{Me^{Kt}},\quad \bar{z}(t)=Me^{Kt}, \quad K>1,\quad M=\max(c_1, c_2),
$$
and verify that the solution of (\ref{FG}) is global existence. Consequently, the solution of (\ref{F-G}) is global existence for initial data $u_0>1$ and $v_0\geq 0$.

(ii) $f(z)=e^{e^z}$, $g(w)=\ln w$. Although we cannot write out the explicit expression of $F^{-1}$ and $G^{-1}$, we can construct sub-solution of (\ref{FG})
having the form of
$$
\underline{w}(t)=e^{\frac{a}{(1-ct)^K}},\quad \underline{z}(t)=\ln \frac{b}{(1-ct)^L},\quad L>K>1,\quad c\leq \min(\frac{a}{L}, \frac{b}{a}),
$$
and prove that the solution of (\ref{FG}) blows up in finite time. Consequently, the solution of (\ref{F-G}) blows up in finite time for large initial data.

(iii) $f(z)=e^{e^z}$, $g(w)=\ln (\ln w)$. Although we cannot write out the explicit expression of $F^{-1}$ and $G^{-1}$, we can construct sup-solution of (\ref{FG}) having the form of
$$
\bar{w}(t)=e^{e^{Me^{Kt}}},\quad \bar{z}(t)=Me^{Kt}, \quad K>1,\quad M=\max(c_1, c_2),
$$
and verify that the solution of (\ref{FG}) is global existence. Consequently, the solution of (\ref{F-G}) is global existence for initial data $u_0>1$ and $v_0\geq 0$.

(iv) Generally, consider the ODE problem
\begin{equation}
\label{SM11'}\left\{
\begin{array}{llll}
&w_t=\exp(\exp(...(\exp(z))...)),\quad z_t=\ln(\ln(...(\ln(M+w))...)),\quad t>0,\\
&w(0)=c_1,\quad z(0)=c_2.
\end{array}\right.
\end{equation}
Here $\exp(\exp(...(\exp(z))...))$ is $m$-multiple contained function, $\ln(\ln(...(\ln(M+w))...))$ is $n$-multiple contained function, $M$ is large enough such that $$\ln(\ln(...(\ln(M))...))\geq 0.$$

If $m>n$, and $c_1$, $c_2$ are large enough, we can construct the blowup sub-solution having the form of $$
\underline{w}(t)=\exp(\exp(...(\exp(\frac{a}{(1-ct)^K}))...)),\  \underline{z}(t)=\ln(\frac{b}{(1-ct)^L}).$$
Here $\exp(\exp(...(\exp(s))...))$ is $(m-1)$-multiple contained function.

If $m\leq n$, we can construct the global sup-solution having the form of
$$
\bar{w}(t)=\exp(\exp(...(M\exp(Kt)...)),\  \bar{z}(t)=Me^{Lt}.
$$
Here $\exp(\exp(...(\exp(s))...))$ is $(m+1)$-multiple contained function.

\subsection{The proofs of Theorem A and Theorem 1}
\qquad In this subsection, we give the proofs of Theorem A and Theorem 1.

{\bf The proof of Theorem A:} (1). We only prove it in the case of $\int^{+\infty}_{c_1}\frac{dw}{g_1(w)}<+\infty$ and $(c_1,c_2)>(0,0)$. The proof in the case of $\int^{+\infty}_{c_2}\frac{dz}{f_2(z)}<+\infty$ is similar.

Since $f_i(z), f'_i(z)\geq 0$ for $z\geq 0$ and $f_i(z), f'_i(z)\not\equiv 0$
in any subinterval of $(0,+\infty)$, $g_i(w), g'_i(w)\geq 0$ for $w\geq 0$ and $g_i(w), g'_i(w)\not\equiv 0$
in any subinterval of $(0,+\infty)$, $i=1,2$, we have $w_t\geq 0$ and $z_t\geq 0$ for $t>0$ and can take $t_0$ small enough such that $0<c_1=w(0)<w(t_0)<+\infty$, $0<c_2=z(0)<z(t_0)<+\infty$, $f(z(t))>f(z(t_0))>f(z(0))\geq 0$ for $t>t_0>0$. (\ref{661}) implies that $f_1(z(t))g_1(w(t))>f_1(z(t_0))g_1(w(t))$ and
\begin{align}
w_t=f_1(z(t))g_1(w(t))>f_1(z(t_0))g_1(w(t)),\quad t>t_0,\quad w(t_0)>c_1>0.\label{10191}
\end{align}
But the solution of $$w_t=f_1(z(t_0))g_1(w(t)),\quad t>t_0,\quad w(t_0)>c_1>0$$
will blow up in finite time because of $\int^{+\infty}_{c_1}\frac{dw}{g_1(w)}<+\infty$. Therefore, the solution of (\ref{661}) will blow up in finite time.

(2). Noticing that (\ref{662})--(\ref{661''}) and (\ref{661''}) has the form of (\ref{FG}), by the results and the proof of Theorem A$'$, the conclusions of Theorem A are true.\hfill $\Box$

{\bf The proof of Theorem 1:} (1). Let $(\bar{u}(t),\bar{v}(t))$ be the global solution of (\ref{1'}) with initial data
$$
c_1=\max_{x\in\bar{\Omega}}u_0(x),\quad c_2=\max_{x\in\bar{\Omega}}v_0(x).
$$
  It is to verify that $(\bar{u}(t),\bar{v}(t))$ is a supsolution of (\ref{1})(or (\ref{h2})) and
$$
u(x,t)\leq \bar{u}(t),\quad v(x,t)\leq \bar{v}(t)\quad {\rm for }\quad x\in \Omega,\ t>0,
$$
which implies that the solutions of (\ref{1}) and (\ref{h2}) are global existence.

(2). Let $(\underline{u}(t),\underline{v}(t))$ be the blowup solution of (\ref{1'}) with initial data
$$
c_1\leq\min_{x\in\bar{\Omega}}u_0(x),\quad c_2\leq \min_{x\in\bar{\Omega}}v_0(x).
$$
 Since $\frac{\partial f}{\partial v}\geq 0$
and $\frac{\partial g}{\partial u}\geq 0$, we can apply the comparison principle to (\ref{1})(or (\ref{h2})). It is easy to verify that $(\underline{u}(t),\underline{v}(t))$ is a sub-solution of (\ref{1}) (or (\ref{h2})) and
$$
u(x,t)\geq \underline{u}(t),\quad v(x,t)\geq \underline{v}(t)\quad {\rm for }\quad x\in \Omega,\ t>0,
$$
which implies that the solutions of (\ref{1}) and (\ref{h2}) will blow up in finite time. \hfill$\Box$

{\bf Example 2.4.} 1. $f_1(v)g_1(u)=v^pu^q$, $f_2(v)g_2(u)=u^rv^s$, $pr=(1-q)(1-s)$ is the watershed for the blowup in finite time and global existence of the solution.

2. $f_1(v)g_1(u)=\frac{v^p}{u^q}$, $f_2(v)g_2(u)=\frac{u^r}{v^s}$, $pr=(q+1)(s+1)$ is the watershed for the blowup in finite time and global existence of the solution.

{\bf Example 2.5.} Consider
 \begin{equation}
\label{SM11}\left\{
\begin{array}{llll}
&u_t=d_1\Delta u+\exp(\exp(...(\exp(v))...)),\quad x\in\Omega,\ t>0,\\
&v_t=d_2\Delta v+\ln(\ln(...(\ln(M+u))...)),\quad x\in\Omega,\ t>0,\\
&\frac{\partial u}{\partial \eta}=\frac{\partial v}{\partial \eta}=0\quad {\rm or}\ u=v=0,\qquad  \qquad x\in\partial\Omega,\ t>0,\\
&u(x,0)=u_0(x),\quad v(x,0)=v_0(x),\quad x\in\Omega.
\end{array}\right.
\end{equation}
Here $\exp(\exp(...(\exp(v))...))$ is $m$-multiple contained function, $\ln(\ln(...(\ln(M+u))...))$ is $n$-multiple contained function, $M$ is large enough such that $$\ln(\ln(...(\ln(M))...))\geq 0.$$

By the results of Example 2.3 (iv) and Theorem 1, the solution of (\ref{SM11}) is global existence for any nonnegative initial data if $n\geq m$, while the solution will blow up in finite time for some initial data if $n<m$, the watershed is $m=n$.

\section{The results on (\ref{1}) when $f$ and $g$ aren't separation of variables}
\qquad In this section, we will consider (\ref{1}) when $f$ and $g$ aren't separation of variables. Due to the complexity of nonlinearities, whether the solution is global existence or blowup cannot be judged directly by the structures of $f$ and $g$, we need some auxiliary functions.
We have

{\bf Theorem B} {\it (1). Assume that $f(w,z)$ and $g(w,z)$ are locally Lipschitz continuous functions for $(w,z)\geq \mathbf{0}$, and there exist nonnegative smooth functions $h(w)$ and $l(z)$ such that
\begin{align}
f(w,z)h'(w)l(z)+g(w,z)h(w)l'(z)\leq H[h(w)l(z)],\label{12293}
\end{align}
 where $H(s)$ is a continuous function for $s\geq 0$. If $H(s)>0$ for $s\geq 0$ and
 \begin{align}
 \int^{+\infty}_0\frac{ds}{H(s)}=+\infty,\label{1229w2}
 \end{align}
or $-\infty<H(s)\leq 0$ for $s\geq 0$, then $h(w)l(z)$ is global existence.

Especially, if (\ref{1229w2}) holds, $f(w,z)\geq 0$, $g(w,z)\geq 0$, $h'(w)\geq 0(\not\equiv 0)$ and $l'(z)\geq 0(\not\equiv 0)$, then the solution of (\ref{1'}) is global existence for any nonnegative initial data.

(2). Assume that $f(w,z)$ and $g(w,z)$ are locally Lipchitz continuous functions for all $(w,z)\geq \mathbf{0}$, and there exist nonnegative smooth functions $h(w)$ and $l(z)$  such that
\begin{align}
f(w,z)h'(w)l(z)+g(w,z)h(w)l'(z)\geq L[h(w)l(z)]\geq 0,\label{12291}
\end{align}
where $L(s)$ is a continuous function for $s\geq 0$. If
 \begin{align}
 \int^{+\infty}_c\frac{ds}{L(s)}<+\infty,\quad c>0,\label{1229w3}
 \end{align}
  then $h(w)l(z)$ and the solution of (\ref{1'}) will blow up in finite time for large initial data.}

{\bf Proof:} (1). Multiplying the equation $w_t=f(w,z)$ by $h'(w)l(z)$, the equation $z_t=g(w,z)$ by $h(w)l'(z)$, then adding the results, we have
\begin{align}
h'(w)l(z)w_t+h(w)l'(z)z_t=f(w,z)h'(w)l(z)+g(w,z)h(w)l'(z)\leq H(h(w)l(z)).\label{1229w1}
\end{align}
By the classic results of ODE, if $H(s)\geq 0$ with the condition of (\ref{1229w2}) or $-\infty<H(s)\leq 0$, then the solution of the following ODE problem
\begin{align}
s_t=H(s(t)),\quad s(0)=c>0\label{12294}
\end{align}
is global existence. However, the solution of ({\ref{12294}) with $s(0)=h(c_1)l(c_2)$ is a sup-solution of (\ref{1229w1}).
 Consequently, $h(w(t))l(z(t))$ is global existence.

 Especially, if (\ref{1229w2}) holds, $f(w,z)$, $g(w,z)$, $h'(w)$ and $l'(z)$ are nonnegative, then we have
 $w_t\geq 0$, $z_t\geq 0$, $h(w)\geq h(w(t_0))>0$ and $l(z)\geq l(z(t_0))>0$ for some $t_0>0$. Therefore,
 $$
 h(w(t_0))l(z(t))\leq h(w(t))l(z(t)),\quad  h(w(t))l(z(t_0))\leq h(w(t))l(z(t)),
 $$
which implies that $h(w(t))$ and $l(z(t))$ are global existence. By the continuities and monotonicity of $h(w)$ and $l(z)$, we know that $(w(t),z(t))$ is global existence for any nonnegative initial data.

(2). Multiplying the equation $w_t=f(w,z)$ by $h'(w)l(z)$, the equation $z_t=g(w,z)$ by $h(w)l'(z)$, then adding the results, we have
\begin{align}
h'(w)l(z)w_t+h(w)l'(z)z_t=f(w,z)h'(w)l(z)+g(w,z)h(w)l'(z)\geq L(h(w)l(z)).\label{1229w4}
\end{align}
By the classic theory of ODE, under the condition of (\ref{1229w3}), the solution of the following ODE problem
\begin{align}
s_t=L(s(t)),\quad s(0)=c>0\label{12295}
\end{align}
will blow up in finite time for large initial data. The solution of  (\ref{12295}) with $s(0)=h(c_1)l(c_2)$ is a sub-solution of (\ref{1229w4}). Consequently,
$h(w(t))l(z(t))$ will blow up in finite time for large initial data. By the continuities of $h(w)$ and $l(z)$, we know that $(w(t),z(t))$  will blow up in finite time for large initial data.\hfill$\Box$

{\bf Example 3.1.} $f(w,z)=Aw^2z-Bw^p$, $g(w,z)=Cwz^2-Dz^q$($A, B, C, D>0$). If $p\neq 2$ and $q\neq 2$, we cannot write $f(w,z)$ and $g(w,z)$ in the form of
$f(w,z)=f_1(z)g_1(w)$ and $g(w,z)=f_2(z)g_2(w)$.

(1). However,  if $0<p<2$ and $0<q<2$, and $(w(0), z(0))>(M,N)$ with $(M,N)$ is large enough such that $f(w(0),z(0))\geq 0$ and $g(w(0),z(0))\geq 0$,  taking $h(w)=w$ and $l(z)=z$, then we can get
\begin{align}
(wz)_t=(A+C)w^2z^2-Bw^pz-Dwz^q,\label{9151}
\end{align}
and prove that $wz$ and the solution of (\ref{1'}) will blow up in finite time for large initial data. In fact, by continuities, we have $w_t\geq 0$ and $z_t\geq 0$ for all $t\geq 0$, and we can write (\ref{9151}) as
\begin{align}
(wz)_t=w^2z^2[(A+C)-\frac{B}{w^{2-p}z}-\frac{D}{wz^{2-q}}]\geq \frac{(A+C)}{2}(wz)^2:=L(wz).
\end{align}
Obviously, $(wz)$ will blow up in finite time for large initial data.
Moreover, we even can construct the sub-solution of (\ref{1'}) has the form of $$(\underline{u}(t),\ \underline{v}(t))=(\frac{M}{(1-ct)},\ \frac{N}{(1-ct)})$$
with $w(0)=M$, $z(0)=N$ and $$AM^2N>BM^p,\quad CMN^2>DN^q$$ if $0<p<2$ and $0<q<2$. Here $$c=\min (AMN-BM^{p-1},\ CMN-DN^{q-1}).$$

(2). On the other hand, if $0<p<2$, $0<q<2$, the initial data $w(0)=\epsilon_1$ and $z(0)=\epsilon_2$ are small enough such that $\epsilon_1^{2-p}\epsilon_2\leq \frac{B}{A}$ and $\epsilon_1\epsilon_2^{2-q}\leq \frac{D}{C}$, then (\ref{9151}) can be written as
\begin{align}
(wz)_t=w^pz(Aw^{2-p}z-B)+wz^q(Cwz^{2-q}-D)\leq 0\equiv H(wz),\label{1171}
\end{align}
the solution of (\ref{1'}) is global existence. In fact, we even can let $(\epsilon_1, \epsilon_2)$ be a supsolution of (\ref{1'}), and obtain $w(t)\leq \epsilon_1$, $z(t)\leq \epsilon_2$.

{\bf Example 3.2.} (1). The case of $f(w,z)\leq 0$ and $g(w,z)\geq 0$ is contained in the cases of Theorem B. For example: Consider the following problem
\begin{equation}
\label{915x1}\left\{
\begin{array}{llll}
&w_t=-Aw^k(t)z^l(t),\quad z_t=Bw^m(t)z^n(t),\quad t>0\\
&w(0)=c_1\geq 0,\quad z(0)=c_2\geq 0.
\end{array}\right.
\end{equation}
If there exist $\alpha>0$ and $\beta>0$ such that
\begin{align}
\frac{k+\alpha-1}{\alpha}=\frac{l+\beta}{\beta}:=p,\quad \frac{m+\alpha}{\alpha}=\frac{n+\beta-1}{\beta}:=q.\label{915x2}
\end{align}
Then
\begin{align}
(w^{\alpha}z^{\beta})_t=-A\alpha (w^{\alpha}z^{\beta})^p+B\beta(w^{\alpha}z^{\beta})^q.\label{915x3}
\end{align}
 If $1<p\leq q$ and $A\alpha<B\beta$, $(w^{\alpha}z^{\beta})$ will blow up in finite time for positive initial data .

(2). In Theorem B and Theorem 2, the conditions $f(w,z)\geq 0$, $g(w,z)\geq 0$, $h'(w)\geq 0$ and $l'(z)\geq 0$ are needed to keep the solution $(w,z)$ of (\ref{1'}) be global existence. We consider (\ref{915x1}) again, but we will show that: There exist $A, B, k, l, m$ and $l$ such that
$\lim_{t\rightarrow T^-} w(t)=0$, $\lim_{t\rightarrow T^-} z(t)=+\infty$ but $wz$ keeps bounded in $[0,T)$.

Let
$$
w(t)=(T-t)^a,\quad z(t)=\frac{1}{(T-t)^b}.
$$
Suppose that $a$, $b$, $k$, $l$, $m$ and $n$ satisfy
$$
a(1-k)+bl=1,\quad am+b(1-n)=-1
$$
and
$$
a=\frac{l+1-n}{(1-k)(1-n)-lm}>0,\quad b=\frac{k-1-m}{(1-k)(1-n)-lm}>0.
$$
(For example, $k=1$, $l=2$, $m=1$, $n=4$, $a=b=\frac{1}{2}$).

If $A=a$, $B=b$ and initial data
$$
c_1=T^a,\quad c_2=\frac{1}{T^b},
$$
then (\ref{915x1}) has the exact solution $(w(t),z(t))=((T-t)^a, \frac{1}{(T-t)^b})$. Moreover, $\lim_{t\rightarrow T^-} w(t)=0$, $\lim_{t\rightarrow T^-} z(t)=+\infty$ but $wz$ keeps bounded in $[0,T)$ if $a\geq b>0$.

Parallel to Theorem B, we have the following conclusions on (\ref{1}).

{\bf Theorem 2.} {\it  Assume that $f(u,v)$ and $g(u,v)$ are locally Lipchitz continuous functions for $(u,v)\geq \mathbf{0}$, $\frac{\partial f(u,v)}{\partial v}\geq 0$ and $\frac{\partial g(u,v)}{\partial u}\geq 0$. Suppose that $h(u)\geq 0$ and $h'(u)\geq 0$ for $u\geq 0$, $l(v)\geq 0$ and $l'(v)\geq 0$ for $v\geq 0$.

(1). If (\ref{12293}) and (\ref{1229w2}) hold, then the solutions of (\ref{1}A), (\ref{1}B), (\ref{h2}A) and (\ref{h2}B) are global existence for any nonnegative initial data;

(2). Suppose that (\ref{12041}) and (\ref{12042}) hold. If (\ref{12291}) and (\ref{1229w3}) hold and the solution of (\ref{1'}) blows up in finite time for initial data $(c_1,c_2)$, then the solutions of (\ref{1}A) and (\ref{h2}A) will blow up in finite time for initial data $(u_0(x),v_0(x))\geq (c_1, c_2)$. If the solution of (\ref{1'}) blows up in finite time for initial data $(c_1,c_2)=\mathbf{0}$, then the solutions of (\ref{1}B) and (\ref{h2}B) will blow up in finite time for initial data $(u_0(x),v_0(x))\geq \bf{0}$.
}

{\bf Proof:} The proof is entirely similar to that of Theorem 1. We can take the solution of (\ref{1'}) as the sub-solution or sup-solution of (\ref{1}) and (\ref{h2}), then we can obtain the corresponding conclusions.\hfill $\Box$

\section{Quenching phenomenon}
\qquad In this section, we discuss the conditions on $f(u,v)$ and $g(u,v)$ which can make the solutions of (\ref{1}), (\ref{1'}) and (\ref{h2}) quench in finite time.

{\bf Theorem C} {\it (1). Assume that there exist $0\leq c'_1\leq c_1$, $0\leq c'_2\leq c_2$ such that $f(w,z)$ and $g(w,z)$ are locally Lipchitz continuous functions,  $f(w,z)< 0$ and $g(w,z)<0$ for $c'_1<w\leq c_1$ and $c'_2<z\leq c_2$. If
\begin{align}
\lim_{w\rightarrow c'_1}|f(w,z)|=+\infty\quad {\rm for}\ z\geq 0,\label{12231}
\end{align}
or
 \begin{align}
\lim_{z\rightarrow c'_2}|g(w,z)|=+\infty\quad {\rm for}\ w\geq 0,\label{12232}
\end{align}
or
\begin{align}
\lim_{(w,z)\rightarrow (c'_1,c'_2)}|f(w,z)|=+\infty,\quad
\lim_{(w,z)\rightarrow (c'_1,c'_2)}|g(w,z)|=+\infty,\label{12233}
\end{align}
then the solution of (\ref{1'}) will quench in finite time.

(2). Assume that there exist $c_1\leq c''_1<+\infty$,  $c_2\leq c''_2<+\infty$ such that $f(w,z)$ and $g(w,z)$ are locally Lipchitz continuous functions,  $f(w,z)>0$ and $g(w,z)>0$ for $c_1<w\leq c''_1$ and $c_2<z\leq c''_2$. If
\begin{align}
\lim_{w\rightarrow c''_1}|f(w,z)|=+\infty\quad {\rm for}\ z\geq 0,\label{66w1}
\end{align}
or
 \begin{align}
\lim_{z\rightarrow c''_2}|g(w,z)|=+\infty\quad {\rm for}\ w\geq 0,\label{66w2}
\end{align}
or
\begin{align}
\lim_{(w,z)\rightarrow (c''_1,c''_2)}|f(w,z)|=+\infty,\quad
\lim_{(w,z)\rightarrow (c''_1,c''_2)}|g(w,z)|=+\infty,\label{66w3}
\end{align}
then the solution of (\ref{1'}) will quench in finite time.
}

{\bf Proof:} (1). Since $f(w,z)<0$ and $g(w,z)<0$, we have $w_t<0$, $z_t<0$, $w(t)\leq c_1$ and $z(t)\leq c_2$ in the time interval $[0,t]$ wherever they exist.

On the other hand, by (\ref{12231})--(\ref{12233}), there exists $T>0$ such that at least one of the following two cases holds:

Case (i) $w(t)\rightarrow c'_1$ and $|w_t|\rightarrow +\infty$ as $t\rightarrow T^-$;

Case (ii) $z(t)\rightarrow c'_2$ and $|z_t|\rightarrow +\infty$ as $t\rightarrow T^-$.

In any case, quenching phenomenon happens.

(2). Since $f(w,z)>0$ and $g(w,z)>0$ for $c_1<w\leq c''_1$ and $c_2<z\leq c''_2$, we have $w_t>0$, $z_t>0$, $w(t)\geq c_1$ and $z(t)\geq c_2$ in the time interval $[0,t]$ wherever they exist.

On the other hand, by (\ref{66w1})--(\ref{66w3}), there exists $T>0$ such that at least one of the following two cases holds:

Case (i) $w(t)\rightarrow c''_1$ and $|w_t|\rightarrow +\infty$ as $t\rightarrow T^-$;

Case (ii) $z(t)\rightarrow c''_2$ and $|z_t|\rightarrow +\infty$ as $t\rightarrow T^-$.

In any case, quenching phenomenon occurs.\hfill $\Box$

Parallel to Theorem C, we have the following conclusions.

{\bf Theorem 3.} {\it (1). Assume that there exist $0\leq c'_1\leq \min_{x\in \Omega}u_0(x)$,  $0\leq c'_2\leq \min_{x\in \Omega}v_0(x)$ such that $f(u,v)$ and $g(u,v)$ are locally Lipchitz continuous functions, $f(u,v)<0$, $g(u,v)<0$, $\frac{\partial f}{\partial v}>0$ and $\frac{\partial g}{\partial u}>0$ for $c'_1<u\leq \max_{x\in \Omega}u_0(x)$ and $c'_2<v\leq \max_{x\in \Omega}v_0(x)$. Moreover, suppose that
\begin{align}
&d_1\Delta u_0-a(x)\cdot \nabla u_0+f(u_0,v_0)\leq 0\quad {\rm for}\quad  x\in\Omega,\label{12045}\\
&d_2\Delta v_0-b(x)\cdot \nabla v_0+g(u_0,v_0)\leq 0\quad {\rm for}\quad  x\in\Omega.\label{12046}
\end{align}
If one of (\ref{12231})--(\ref{12233}) holds, then the solutions of (\ref{1}) and (\ref{h2}) will quench in finite time. Moreover, if $f(u,v)\leq -a_1-a_2u$ and there exists $A(x)$ such that $\varphi_1(x)a(x)=\nabla A(x)$ and $\Delta A(x)\leq -a_3$, where $a_1, a_2>0$, $a_3\geq 0$, then the quenching time satisfies
\begin{align}
T_{\max}&\leq
\frac{1}{\lambda_1d_1+a_2}\ln \frac{a_3c_1|\Omega|+a_1+(\lambda_1d_1+a_2)\int_{\Omega} u_0(x)\varphi_1(x)dx}{a_3c_1|\Omega|+a_1+(\lambda_1d_1+a_2)c'_1}.\label{T1}
\end{align}
If $g(u,v)\leq -b_1-b_2v$ and there exists $B(x)$ such that $\varphi_1(x)b(x)=\nabla B(x)$ and $\Delta B(x)\leq -b_3$, where $b_1, b_2>0$, $b_3\geq 0$, then the quenching time satisfies
\begin{align}
T_{\max}&\leq
\frac{1}{\lambda_1d_2+b_2}\ln \frac{b_3c_2|\Omega|+b_1+(\lambda_1d_2+b_2)\int_{\Omega} v_0(x)\varphi_1(x)dx}{b_3c_2|\Omega|+b_1+(\lambda_1d_2+b_2)c'_2}.\label{T1'}
\end{align}
Here $\varphi_1$ is the first eigenfunction(normalized by $\int_{\Omega} \varphi_1(x)dx=1$) of the following eigenvalue problem:
\begin{align}
-\Delta \varphi=\lambda_1 \varphi\quad {\rm in}\ \Omega,\quad \varphi=0\quad {\rm on}\quad \partial\Omega. \label{tzz}
\end{align}

(2). Assume that there exist $\max_{x\in \Omega}u_0(x)\leq c''_1<+\infty$,  $\max_{x\in \Omega}v_0(x)\leq c''_2<+\infty$ such that $f(u,v)$ and $g(u,v)$ are locally Lipchitz continuous functions, $f(u,v)>0$, $g(u,v)>0$, $\frac{\partial f}{\partial v}>0$ and $\frac{\partial g}{\partial u}>0$ for $0\leq u\leq c''_1$ and $0\leq v\leq c''_2$. Moreover, suppose that (\ref{12041}) and (\ref{12042}) hold. If one of (\ref{66w1})--(\ref{66w3}) holds, then the solutions of (\ref{1}) and (\ref{h2}) will quench in finite time. Moreover, if $f(u,v)\geq a'_1+a'_2u$, $a'_2>\lambda_1d_1$, there exists $A(x)$ such that $\varphi_1(x)a(x)=\nabla A(x)$ and $\Delta A(x)\geq a'_3\geq 0$, then the quenching time satisfies
\begin{align}
T_{\max}&\leq
\frac{1}{(a'_2-\lambda_1d_1)}\ln \frac{a'_3c''_1|\Omega|+a'_1+(a'_2-\lambda_1d_1)c''_1}{a'_3c''_1|\Omega|+a'_1+(a'_2-\lambda_1d_1)\int_{\Omega} v_0(x)\varphi_1(x)dx}.\label{T2}
\end{align}
If $g(u,v)\geq b'_1+b'_2v$, $b'_2>\lambda_1d_2$, there exists $B(x)$ such that $\varphi_1(x)b(x)=\nabla B(x)$ and $\Delta B(x)\geq b'_3\geq 0$, then
the quenching time satisfies
\begin{align}
T_{\max}&\leq
\frac{1}{(b'_2-\lambda_1d_2)}\ln \frac{b'_3c''_2|\Omega|+b'_1+(b'_2-\lambda_1d_2)c''_2}{b'_3c''_2|\Omega|+b'_1+(b'_2-\lambda_1d_2)\int_{\Omega} v_0(x)\varphi_1(x)dx}.\label{T2'}
\end{align}
}

{\bf Proof:} (1). Without loss of generality, we assume that (\ref{12231}) holds.
Since $\frac{\partial f}{\partial v}>0$ and $\frac{\partial g}{\partial u}>0$, the comparison principle can be applied to (\ref{1})(or(\ref{h2})), we can take the solution of (\ref{1'}) with $c_1=\max_{x\in \Omega} u_0(x)$, $c_2=\max_{x\in \Omega} v_0(x)$ as a sup-solution of (\ref{1})(or(\ref{h2})). (\ref{12045}) and (\ref{12046}) imply that there exists $T_{\max}<T$($T$ is the quenching time of (\ref{1'})) such that
\begin{align}
\lim_{t\rightarrow T^-_{\max}} \inf_{x\in \Omega} u(x,t)=c_1.\label{9141}
\end{align}
Consequently,
$$
\lim_{t\rightarrow T^-_{\max}}\sup_{x\in \Omega}|u_t-d_1\Delta u+a(x)\cdot \nabla u|=\lim_{t\rightarrow T^-_{\max}}\sup_{x\in \Omega}|f(u(x,t),v(x,t))|=+\infty,
$$
and at least one of the following equalities holds:

(i) $\lim_{t\rightarrow T^-_{\max}}\sup_{x\in \Omega}|u_t(x,t)|=+\infty$; (ii) $\lim_{t\rightarrow T^-_{\max}}\sup_{x\in \Omega}|\Delta u(x,t)|=+\infty$;
(iii) $\lim_{t\rightarrow T^-_{\max}}\sup_{x\in \Omega}|\nabla u(x,t)|=+\infty$.

 That is, quenching phenomenon occurs.

To give the estimate for the quenching time, multiplying the first equation of (\ref{1})(or (\ref{h2})) by $\varphi_1(x)$ and integrating it on $\Omega$, we get
\begin{align}
\frac{d}{dt}\int_{\Omega} u\varphi_1(x)dx&=-\lambda_1d_1\int_{\Omega} u\varphi_1(x)dx-\int_{\Omega}\varphi_1(x) a(x)\cdot \nabla udx+\int_{\Omega}f(u,v)\varphi_1(x)dx\nonumber\\
&\leq -\lambda_1d_1\int_{\Omega} u\varphi_1(x)dx-\int_{\Omega}\nabla A(x)\cdot \nabla udx+\int_{\Omega}(-a_1-a_2u)\varphi_1(x)dx\nonumber\\
&= -(\lambda_1d_1+a_2)\int_{\Omega} u\varphi_1(x)dx+\int_{\Omega}\Delta A(x) udx-\int_{\Omega}a_1\varphi_1(x)dx\nonumber\\
&\leq -(\lambda_1d_1+a_2)\int_{\Omega} u\varphi_1(x)dx-a_3c_1|\Omega|-a_1.\label{611}
\end{align}
Letting
$$
y(t)=\int_{\Omega} u\varphi_1(x)dx,
$$
using (\ref{611}), we have
$$
\frac{dy}{dt}\leq -(\lambda_1d_1+a_2)y(t)-a_3c_1|\Omega|-a_1, \quad {\rm i.e.} \quad \frac{dy}{(\lambda_1d_1+a_2)y(t)+a_3c_1|\Omega|+a_1}\leq -dt.
$$
Integrating it on $[0,T_{\max})$, we get
\begin{align*}
T_{\max}&\leq \frac{1}{\lambda_1d_1+a_2}\ln \frac{a_3c_1|\Omega|+a_1+(\lambda_1d_1+a_2)y(0)}{a_3c_1|\Omega|+a_1+(\lambda_1d_1+a_2)y(T_{\max})}\\
&\leq
\frac{1}{\lambda_1d_1+a_2}\ln \frac{a_3c_1|\Omega|+a_1+(\lambda_1d_1+a_2)\int_{\Omega} u_0(x)\varphi_1(x)dx}{a_3c_1|\Omega|+a_1+(\lambda_1d_1+a_2)c'_1}.
\end{align*}

(2). Similarly, without loss of generality, we assume that (\ref{66w2}) holds.
Since $\frac{\partial f}{\partial v}>0$ and $\frac{\partial g}{\partial u}>0$, the comparison principle can be applied to (\ref{1})(or(\ref{h2})), we can take the solution of (\ref{1'}) with $c_1=\min_{x\in \Omega} u_0(x)$, $c_2=\min_{x\in \Omega} v_0(x)$ as a sub-solution of (\ref{1})(or(\ref{h2})). (\ref{12041}) and (\ref{12042}) imply that there exists $T_{\max}<T$($T$ is the quenching time of (\ref{1'})) such that
\begin{align}
\lim_{t\rightarrow T^-_{\max}} \sup_{x\in \Omega} v(x,t)=c''_2.\label{9141'}
\end{align}
Consequently,
$$
\lim_{t\rightarrow T^-_{\max}}\sup_{x\in \Omega}|v_t-d_2\Delta v+b(x)\cdot \nabla v|=\lim_{t\rightarrow T^-_{\max}}|g(\sup_{x\in \Omega}(u(x,t),v(x,t)))|=+\infty,
$$
and at least one of the following equalities holds:

(i) $\lim_{t\rightarrow T^-_{\max}}\sup_{x\in \Omega}|v_t(x,t)|=+\infty$; (ii) $\lim_{t\rightarrow T^-_{\max}}\sup_{x\in \Omega}|\Delta v(x,t)|=+\infty$;
(iii) $\lim_{t\rightarrow T^-_{\max}}\sup_{x\in \Omega}|\nabla v(x,t)|=+\infty$.

 That is, quenching phenomenon happens.

Now we give the estimate for the quenching time below. Multiplying the second equation of (\ref{1})(or (\ref{h2})) by $\varphi_1(x)$ and integrating it on $\Omega$, we have
\begin{align}
\frac{d}{dt}\int_{\Omega} v\varphi_1(x)dx&=-\lambda_1d_2\int_{\Omega} v\varphi_1(x)dx-\int_{\Omega}\varphi_1(x) b(x)\cdot \nabla vdx+\int_{\Omega}g(u,v)\varphi_1(x)dx\nonumber\\
&\geq -\lambda_1d_2\int_{\Omega} v\varphi_1(x)dx-\int_{\Omega}\nabla B(x)\cdot \nabla vdx+\int_{\Omega}(b'_1+b'_2v)\varphi_1(x)dx\nonumber\\
&= (b'_2-\lambda_1d_2)\int_{\Omega} v\varphi_1(x)dx+\int_{\Omega}\Delta B(x) vdx+\int_{\Omega}b'_1\varphi_1(x)dx\nonumber\\
&\geq (b'_2-\lambda_1d_2)\int_{\Omega} v\varphi_1(x)dx+b'_3c''_2|\Omega|+b'_1.\label{9142'}
\end{align}
Letting
$$
\tilde{y}(t)=\int_{\Omega} v\varphi_1(x)dx,
$$
using (\ref{9142'}), we get
$$
\frac{d\tilde{y}}{dt}\geq (b'_2-\lambda_1d_2)\tilde{y}+b'_3c''_2|\Omega|+b'_1, \quad {\rm i.e.} \ \frac{d\tilde{y}}{(b'_2-\lambda_1d_2)\tilde{y}+b'_3c''_2|\Omega|+b'_1}\geq dt.
$$
Integrating it on $[0,T_{\max})$, we obtain
\begin{align*}
T_{\max}&\leq \frac{1}{(b'_2-\lambda_1d_2)}\ln \frac{b'_3c''_2|\Omega|+b'_1+(b'_2-\lambda_1d_2)\tilde{y}(T_{\max})}{b'_3c''_2|\Omega|+b'_1+(b'_2-\lambda_1d_2)\tilde{y}(0)}\\
&\leq
\frac{1}{(b'_2-\lambda_1d_2)}\ln \frac{b'_3c''_2|\Omega|+b'_1+(b'_2-\lambda_1d_2)c''_2\int_{\Omega} \varphi_1(x)dx}{b'_3c''_2|\Omega|+b'_1+(b'_2-\lambda_1d_2)\int_{\Omega} v_0(x)\varphi_1(x)dx}.
\end{align*}

Theorem 3 is proved. \hfill$\Box$

{\bf Example 4.1.} 1. $f(u,v)=-u^{-p_1}v^{-q_1}$, $g(u,v)=-u^{-p_2}v^{-q_2}$, where $p_1\geq 0$, $q_1>0$, $p_2>0$, $q_2\geq 0$, and the initial data $(c_1,c_2)>(0,0)$, then there exists some $T>0$ such that quenching phenomenon happens at $T$.

2. $f(u,v)=\frac{u}{4-v^2}$, $g(u,v)=\frac{v}{2-u}$, and the initial data $(2,2)>(c_1,c_2)>(0,0)$, then there exists some $T>0$ such that quenching phenomenon happens at $T$.

{\bf Remark 4.1.} Since $c'_1<\max_{x\in \Omega}u_0(x)$ and $\int_{\Omega} \varphi_1(x)dx=1$, we have
\begin{align*}
&\quad\frac{1}{\lambda_1d_1+a_2}\ln \frac{a_3c_1|\Omega|+a_1+(\lambda_1d_1+a_2)\int_{\Omega} u_0(x)\varphi_1(x)dx}{a_3c_1|\Omega|+a_1+(\lambda_1d_1+a_2)c'_1}\\
&\leq \frac{1}{\lambda_1d_1+a_2}\ln \frac{a_1+(\lambda_1d_1+a_2)\int_{\Omega} u_0(x)\varphi_1(x)dx}{a_1+(\lambda_1d_1+a_2)c'_1},
\end{align*}
which implies that convection term can delay the quenching time.

\section{Some type convection terms make the solution exist globally}
\qquad In this section, we will show that some type convection terms $a(x)\cdot \nabla u$ and $b(x)\cdot \nabla v$  make the solution exist globally.

{\bf Theorem 4.} {\it  Assume that $f(u,v)$ and $g(u,v)$ are smooth functions for $(u,v)\geq \mathbf{0}$, $\frac{\partial f}{\partial v}\geq 0$
and $\frac{\partial g}{\partial u}\geq 0$.

(1). Suppose that the initial data $(u_0, v_0) \in C^2(\Omega)\times C^2(\Omega)$, $(u_0(x),v_0(x))\geq \mathbf{0}$, $\nabla u_0(x)\neq \mathbf{0}$ and $\nabla v_0(x)\neq \mathbf{0}$ for all $x\in \Omega$. Moreover, for (\ref{1}A) and (\ref{h2}A), suppose that $\frac{\partial u_0}{\partial \eta}\geq 0$ and $\frac{\partial v_0}{\partial \eta}\geq 0$ on $\partial\Omega$. Then there exist $a(x)$ and $b(x)$ such that the solutions of (\ref{1}A) and (\ref{1}B) are global existence.

(2). Assume that there exist positive constants $K>1$, $L>1$, $c_3$, $c_4$, $p_1$, $p_2$, $q_1$ and $q_2$ such that $0<p_1<1$, $0<q_2<1$, $p_2q_1\geq (1-p_1)(1-q_2)$,
$c_3^{p_1}c_4^{q_1}\leq \lambda_1d_1c_3K$, $c_3^{p_2}c_4^{q_2}\leq \lambda_1d_2c_4L$,
\begin{align}
f(c_3\varphi^K,c_4\varphi^L)\leq c_3^{p_1}c_4^{q_1}\varphi^{Kp_1+Lq_1},\quad g(c_3\varphi^K,c_4\varphi^L)\leq c_3^{p_2}c_4^{q_2}\varphi^{Kp_2+Lq_2},\label{12236}
\end{align}
and
\begin{align}
\varphi(x) a(x)\cdot \nabla \varphi(x)\geq (K-1)|\nabla \varphi(x)|^2,\\
 \varphi(x)b(x)\cdot \nabla \varphi(x)\geq (L-1)|\nabla \varphi(x)|^2.
\end{align}
Here $\varphi$ is the first eigenfunction of (\ref{tzz}) normalised by $\max_{x\in \bar{\Omega}} \varphi(x)=1$. Then the solutions of (\ref{1}A) and (\ref{1}B), are global existence if the initial data $(u_0(x),v_0(x))$ satisfies $u_0(x)\leq c_3\varphi^K(x)$ and $v_0(x)\leq c_4\varphi^L(x)$.
}

{\bf Proof:} (1). Note that $\frac{\partial f}{\partial v}\geq 0$
and $\frac{\partial g}{\partial u}\geq 0$, the comparison principle of a system of parabolic equations can be applied to (\ref{1}).
Since $\nabla u_0\neq 0$ and $\nabla v_0\neq 0$ for all $x\in \Omega$, there exist $a(x)$ and $b(x)$ such that
\begin{align}
a(x)\cdot \nabla u_0\geq d_1\Delta u_0+f(u_0,v_0),\quad b(x)\cdot \nabla v_0\geq d_2\Delta v_0+g(u_0,v_0). \label{9131}
\end{align}
We write it as
\begin{align}
0\geq d_1\Delta u_0-a(x)\cdot \nabla u_0+f(u_0,v_0),\quad 0\geq d_2\Delta v_0-b(x)\cdot \nabla v_0+g(u_0,v_0).\label{10251}
\end{align}
Combining this and the assumptions of $\frac{\partial u_0}{\partial \eta}\geq 0$ and $\frac{\partial v_0}{\partial \eta}\geq 0$ on $\partial \Omega$, we know that $(u_0,v_0)$ is a sup-solution of (\ref{1}A) with $a(x)$ and $b(x)$ satisfying (\ref{9131}), hence the solutions of (\ref{1}A) is global existence.

Since $a(x)$ and $b(x)$ satisfy (\ref{10251}), and noticing that $u_0\geq 0$ and $ v_0\geq 0$ on $\partial \Omega$, we can take $(u_0,v_0)$ as a sup-solution of (\ref{1}B). Consequently, the solution of (\ref{1}B) is also global existence(globally bounded).

(2). Since $\frac{\partial f}{\partial v}\geq 0$
and $\frac{\partial g}{\partial u}\geq 0$, we can apply the comparison principle to (\ref{1}B).

Under the assumptions of $f(u,v)$, $g(u,v)$, $a(x)$ and $b(x)$, taking $$(\bar{u}(x,t),\bar{v}(x,t))=(c_3\varphi^K(x), c_4\varphi^L(x)),$$ we have
\begin{align}
&\quad d_1\Delta \bar{u}-a(x)\cdot \nabla \bar{u}+f(\bar{u},\bar{v})\nonumber\\
&=-\lambda_1 d_1c_3 K\varphi^K+c_3K(K-1)\varphi^{K-2}|\nabla \varphi|^2-c_3K\varphi^{K-1}a(x)\cdot \nabla \varphi+f(c_3\varphi^K,c_4\varphi^L)\nonumber\\
&\leq -\lambda_1 d_1c_3 K\varphi^K+c_3K\varphi^{K-2}[(K-1)|\nabla \varphi|^2-\varphi a(x)\cdot \nabla \varphi]+c_3^{p_1}c_4^{q_1}\varphi^{Kp_1+Lq_1}\nonumber\\
&\leq [c_3^{p_1}c_4^{q_1}-\lambda_1d_1 c_3 K]\varphi^K+c_3K\varphi^{K-2}[(K-1)|\nabla \varphi|^2-\varphi a(x)\cdot \nabla \varphi]\nonumber\\
&\leq 0 \quad {\rm for}\quad x\in\Omega,\ t>0,\label{9141}
\end{align}
and
\begin{align}
&\quad d_2\Delta \bar{v}-b(x)\cdot \nabla \bar{v}+g(\bar{u},\bar{v})\nonumber\\
&=-\lambda_1 d_2c_4 L\varphi^L+c_4L(L-1)\varphi^{L-2}|\nabla \varphi|^2-c_4L\varphi^{L-1}b(x)\cdot \nabla \varphi+g(c_3\varphi^K,c_4\varphi^L)\nonumber\\
&\leq -\lambda_1 d_2c_4 L\varphi^L+c_4L\varphi^{L-2}[(L-1)|\nabla \varphi|^2-\varphi b(x)\cdot \nabla \varphi]+c_3^{p_2}c_4^{q_2}\varphi^{Kp_2+Lq_2}\nonumber\\
&\leq [c_3^{p_2}c_4^{q_2}-\lambda_1 d_2c_4 L]\varphi^L+c_4L\varphi^{L-2}[(L-1)|\nabla \varphi|^2-\varphi b(x)\cdot \nabla \varphi]\nonumber\\
&\leq 0\quad {\rm for}\quad x\in\Omega,\ t>0. \label{9142}
\end{align}
Obviously,
\begin{align}
&(\bar{u}(x,t),\bar{v}(x,t))=(c_3\varphi^K(x), c_4\varphi^L(x))=(0,0),\label{9143}\\
&\frac{\partial \bar{u}}{\partial \eta}=c_3K\varphi^{K-1}\frac{\partial \varphi}{\partial \eta}=0,\quad
\frac{\partial \bar{v}}{\partial \eta}=c_4L\varphi^{L-1}\frac{\partial \varphi}{\partial \eta}=0\label{9144}
\end{align}
for $x\in \partial \Omega$ and $t>0$.

(\ref{9141}), (\ref{9142}), (\ref{9143})(or (\ref{9144})) show that $(\bar{u},\bar{v})$ is a sup-solution of (\ref{1}) for initial data
$u_0(x)\leq c_3\varphi^K(x)$ and $v_0(x)\leq c_4\varphi^L(x)$, which means that the solution of (\ref{1}) is globally bounded.\hfill$\Box$

{\bf Remark 5.1.}  Theorem 4 shows that, the convection terms $a(x)\cdot \nabla u$ and $b(x)\cdot \nabla v$ can effect the properties for the solutions under certain conditions.

{\bf Remark 5.2.} The typical $f(u,v)$ and $g(u,v)$ satisfying Theorem 4 are $f(u,v)=u^{p_1}v^{q_1}$ and $g(u,v)=u^{p_2}v^{q_2}$.

\end{document}